\documentclass[11pt]{article}
\usepackage{amsmath,amssymb,mathrsfs}
\usepackage{amsfonts}
\input xy
\xyoption{all}

\newcommand{\pf}{\noindent{\em Proof.}\ }
\newcommand{\qed}{\hfill $\Box$\\}

\newcommand{\ZZ}{\mathbb{Z}}
\newcommand{\QQ}{\mathbb{Q}}

\newcommand{\CC}{\mathbb{C}}
\newcommand{\FF}{\mathbb{F}}
\newcommand{\GG}{\mathbb{G}}

\newcommand{\PP}{\mathbb{P}}
\newcommand{\tensor}{\otimes}

\newcommand{\rc}{\subset}

\newcommand{\isom}{\simeq}

\newcommand{\Gal}{\operatorname{Gal}}

\newcommand{\NS}{\operatorname{NS}}
\newcommand{\Hdg}{\operatorname{Hdg}}
\newcommand{\End}{\operatorname{End}}

\newcommand{\hBr}{\hat{\rm Br}}

\newcommand{\lam}{\lambda}

\begin{document}

\title{K3 Surfaces of Finite Height over Finite Fields \footnote{
	This work was partially supported by
	the Natural Sciences and Engineering Research Council
	of Canada (NSERC).}}
\author{{\sc J.-D. Yu} and {\sc N. Yui}\footnote{
	E-mail addresses: {\tt jdyu@mast.queensu.ca}
				and {\tt yui@mast.queensu.ca}}}
\maketitle

\begin{abstract}
Arithmetic of K3 surfaces defined over finite fields is investigated.
In particular,
we show that any K3 surface $X$ of finite height
over a finite field $k$ of characteristic $p \geq 5$
has a quasi-canonical lifting $Z$ to characteristic $0$,
and that for any such $Z$,
the endormorphism algebra of the transcendental cycles $V(Z)$,
as a Hodge module, is a CM field over $\QQ$.
The Tate conjecture for the product of certain two K3 surfaces is also proved.
We illustrate by examples how to determine explicitly
the formal Brauer group associated to a K3 surface over $k$.
Examples discussed here are all of hypergeometric type.
\end{abstract}

\section{Introduction}

We first recall the Tate conjecture for a smooth projective variety $X$
defined over a finite field.\medskip

{\bf The Tate Conjecture (over finite fields)}
(\cite{T}, Conjecture 1){\bf .}
{\em Let $k$ be a finite field of characteristic $p$, and let $\bar{k}$
be an algebraic closure of $k$.  Let $\Gal(\bar{k}/k)$ be the Galois group
of $\bar{k}$ over $k$.
Let $X$ be an absolutely irreducible smooth projective variety
of dimension $d$ defined over $k$,
and let $X_{\bar{k}}:=X\times_k\bar{k}$.
Then for any prime $\ell\neq p$ and for any integer
$0 \leq r \leq d$,
the $\Gal(\bar{k}/k)$-invariant part
\[ H^{2r}_{et}(X_{\bar{k}}, \QQ_{\ell}(r))^{\Gal(\bar{k}/k)} \]
of the \'etale cohomology group
$H^{2r}_{et}(X_{\bar{k}}, \QQ_{\ell}(r))$
is generated by algebraic cycles
of codimension $r$ in $X$ defined over $k$.}\medskip

Let $k$ be a finite field of characteristic $p$.
Let $X$ be a K3 surface of height $h$ over $k$.
It is known that the Tate conjecture holds true for $X$ over $k$
in the case when $h = 1$ (by Nygaard \cite{N}),
or when $h < \infty$ and $p \geq 5$
(by Nygaard and Ogus \cite{NO}; see Remark 2.3).
The proofs in both cases
depend on constructing a good lifting of $X$
to a complex projective K3 surface $Z$,
called a ``quasi-canonical'' lifting,
such that
the Tate classes in
$H^2_{et}(X_{\bar{k}}, \QQ_{\ell}(1))$
are generated by the Hodge classes
in $H^2(Z, \QQ(1))$
via the comparison
\[ H^2(Z, \QQ(1)) \tensor_{\QQ} \QQ_{\ell}
	= H^2_{et}(X_{\bar{k}}, \QQ_{\ell}(1)). \]
Then the validity of the Tate conjecture follows
because the Hodge classes in $H^2(Z, \QQ(1))$
are generated by divisor classes of $Z$,
and the specialization map to the divisor class group of
$X$ (over $\bar{k}$) provides enough algebraic cycles.\medskip

In this paper, we will study basic properties of quasi-canonical liftings
of K3 surfaces over $k$ of finite height, by making use of the work of Zarhin (\cite{Zarhin_Hodge})
on the Hodge structures for complex K3 surfaces.
Let $X$ be a K3 surface over $k$ of finite height $h$,
and let $Z$ be a quasi-canonical lifting of $X$ to characteristic zero.
We will prove in Section 2 that the endomorphism algebra
of the Hodge structure on the transcendental cycles, $V(Z)$, of $Z$
is a CM field over $\QQ$.
By studying the Frobenius action on various cohomology groups,
we will prove the Tate conjecture for the product of two K3 surfaces
in certain cases. This is done in Section 3.\medskip

We also study the formal Brauer group
$\hBr(X)$
associated to a K3 surface $X$ (defined in \cite{AM})
from the viewpoint of formal group laws.
We compute the formal group laws
explicitly for several examples
including one-parameter families of K3 surfaces.
For this, we follow the method of Stienstra (\cite{St}).
Examples discussed here are all of hypergeometric type in the
sense that the logarithms of  (certain liftings of) $\hBr(X)$
are given by hypergeometric series.
As an application, we study congruences and $p$-adic properties
of the coefficients of the logarithms. Our results on formal group laws
are stated in Section 4. \medskip

{\bf Definitions, notations and conventions.}
If $k$ is a field, $\bar{k}$ denotes
an algebraic closure of $k$.
For a prime $p$,
denote $\ZZ_p$ the $p$-adic integers
and $\QQ_p = \ZZ_p \tensor_{\ZZ} \QQ$.
Throughout the paper,
$\FF_q$ denotes a fixed finite field of $q$ elements.
We also fix an algebraic closure of $\FF_q$ implicitly.\medskip

Suppose $\sigma$ is a linear endomorphism
on a finite dimensional vector space $V$ over a field $k$.
The {\em reciprocal characteristic polynomial of $\sigma$ on $V$}
is the polynomial
\[ \det\left( 1-T\sigma | V \right) \]
with indeterminate $T$. It has coefficients in $k$.
\medskip

Let $k$ be a field.
A smooth, geometrically irreducible,
projective variety $X$ of dimension 2
over $k$ is called a {\em K3 surface}
if the canonical line bundle of $X$
is isomorphic to the structure sheaf $\mathcal{O}_X$
and $H^1(X, \mathcal{O}_X) = 0$.\medskip

Let $X$ be a K3 surface defined over a field $k$.  We fix the following notations.

$\bullet$
$\NS(X) =$ the N\'eron-Severi group of $X$ over $k =$
the classes of divisors over $k$ on $X$
modulo rational equivalence.

$\bullet$
$\rho(X) =$ the rank of $\NS(X)$
as an abelian group.

$\bullet$
Let $k'$ be a field extension of $k$ and
$R$ be a local ring with residue field $k'$.
A {\em lifting $Y$ over $R$ of $X$}
is a flat scheme $Y$ over $R$
such that
$Y \times_R k' = X \times_k k'$.

$\bullet$
Suppose $k$ is a perfect field of characteristic $p > 0$.
Let $\hBr(X)$ denote the formal Brauer
group of $X$. It is known
(see \cite{AM}, Examples in p.90)
that $\hBr(X)$
is representable by a one-dimensional commutative formal group.
Recall that
the {\it height} of $X$ is defined to be the height of $\hBr(X)$
as a formal group.
Let $h$ denote the height of $X$.
Then $h$ can take any integer between $1$ and $10$, or $\infty$
(loc.cit.).
In the former case, $X$ is said to have {\it finite height}.
In particular, if $h=1$, $X$ is said to be {\it ordinary}.
$X$ is called {\em supersingular} if $h = \infty$
(i.e., $\hBr(X)$ is the additive formal group).

\section{A quasi-canonical lifting to characteristic zero}

We fix a finite field $k=\FF_q$ of characteristic $p$, and
also fix a K3 surface $X$ defined over
$k$.
The following definition is a reformulated version
(for a K3 surface)
of the definition of quasi-canonical liftings
given in \cite{NO}, Definition (1.5).
The formulation given here is customized
in more suitable ways for our discussion on K3 surfaces.\medskip

{\bf Definition 2.1.}
A projective K3 surface $Z$ over $\CC$
is called a {\em quasi-canonical lifting}
of $X$ if there exist
\renewcommand{\labelenumi}{{\rm (\alph{enumi})}}
\begin{enumerate}
\item
a complete discrete valuation ring $R \rc \CC$
with residue field $k=\FF_q$,
\item
a K3 surface $Y$ defined over $R$
with $Y_{\CC} = Z$, which lifts $X$, and
\item
an endomorphism $\sigma$ on $H^2(Z, \QQ)$
respecting the Hodge structure
(i.e., a Hodge cycle
in $\End_{\QQ} H^2(Z, \QQ)$),
\end{enumerate}
\smallskip

\noindent satisfying the following two conditions:
\renewcommand{\labelenumi}{{\rm (\arabic{enumi})}}
\begin{enumerate}
\item
Under the natural identifications
\[ H^2(Z, \QQ) \tensor_{\QQ} \QQ_{\ell}
	= H^2(Y_{\CC}, \QQ_{\ell})
	= H^2_{et} (X_{\bar{k}}, \QQ_{\ell}), \]
the endomorphism $\sigma$ coincides with the geometric Frobenius
on $H^2_{et} (X_{\bar{k}}, \QQ_{\ell})$ for any prime $\ell \neq p$.
Similarly,
under the natural identifications
(\cite{F}, 6.1.4)
\[ H^2(Z, \QQ) \tensor_{\QQ} \CC
	= H^2_{dR} (Y/R) \tensor_R \CC
	= H^2_{cris} (X/W) \tensor_W \CC,
	\]
the endormorphism $\sigma$ coincides with
the Frobenius endomorphism on the
crystalline cohomology $H^2_{cris}(X/W)$.
Here $W$ is the ring of Witt vectors of $k$ with the embedding
$W \to R \to \CC$.
\item
There exists an integer $n$ such that the fixed part of $(\sigma/q)^n$
on $H^2(Z, \QQ)$ is precisely equal to
the rational N\'eron--Severi group
$\NS(Z)_{\QQ} := \NS(Z) \tensor_{\ZZ} \QQ$,
and the specialization map
\[ \NS(Z)_{\QQ} = \NS(Y_{\CC})_{\QQ}
	\to \NS(X_{\bar{k}})_{\QQ} \]
is an isomorphism.
\end{enumerate}

We simply call the Hodge cycle $\sigma$ in (c)
the {\em lifted geometric Frobenius} on $Z$.\medskip

Notice that when condition (1) above is fulfilled,
condition (2) implies the validity
of the Tate conjecture
for $X$ over $k$.\medskip

The Tate conjecture for a K3 surface of finite height was
established by Nygaard and Ogus \cite{NO}.\medskip

{\bf Theorem 2.2} (Nygaard and Ogus \cite{NO}){\bf .}
{\em Let $k =\FF_q$ be a finite field of characteristic $p \geq 5$.
Let $X$ be a K3 surface of finite height over $k$.
Then there exist a totally ramified finite extension $R$
of the Witt ring $W(k)$ of $k$
and a lifting $Y$ over $R$ of $X$ such that for any
embedding $R \to \CC$, the induced complex K3 surface $Z = Y_{\CC}$
is a quasi-canonical lifting of $X$.
In particular, the Tate conjecture is true for $X$ over $k$.}\medskip

\pf
Write $W$ for $W(k)$.
By \cite{NO}, Theorem (5.6),
there exist a totally ramified finite extension $R$
of $W$ and a lifting $Y$ over $R$ of $X$
such that the geometric Frobenius $\sigma$
on $H^2_{cris}(X/W)$
preserves the Hodge filtration induced from $Y$
through the identification
\[ H^2_{cris}(X/W) \tensor_W L = H^2_{dR} (Y/R) \tensor_R L. \]
Here $L$ is the field of fractions of $R$.
Then by op.cit., Proposition (1.7),
for any embedding $R \to \CC$,
the map $\sigma$ preserves the rational structure
$H^2(Z, \QQ)$
via the identification
\[ H^2(Z, \QQ) \tensor_{\QQ} \CC = H^2_{dR}(Y/R) \tensor_R \CC. \]
Moreover, this induced map $\sigma$ preserves the Hodge structure
on $H^2(Z, \QQ)$.
By the construction,
$\sigma$ on $H^2(Z, \QQ)$
satisfies condition (1) in Definition 2.1.
Furthermore,
$\sigma$ also satisfies condition (2) in Definition 2.1
(op.cit., Theorem (2.1) and Remark (2.2.3)).
This completes the proof.
\qed

{\bf Remark 2.3.}
The construction of a quasi-canonical lifting of $X$
in \cite{NO}
rests on a fact that there is an equivalence between the category
of the deformations of a K3 surface of finite height,
and that of the deformations of the associated crystal
(see \cite{NO}, Theorem (4.5)).
The assumption that $p \geq 5$
is needed to establish that
the functor between these two categories is fully faithful
(op.cit., Lemma (4.6)).
At this moment, we do not know
if there always exists a quasi-canonical lifting of $X$
when $p = 2$ or 3.\medskip

{\bf (2.4)}
Here we recall some results
from Zarhin \cite{Zarhin_Hodge}.
For any complex projective K3 surface $Z$,
we define the {\em transcendental cycles} of $Z$
to be the orthogonal complement of $\NS(Z)$ in $H^2(Z,\QQ(1))$ with respect to
the cup-product:
$$V(Z):= \NS(Z)_{\QQ}^{\perp} \rc H^2(Z, \QQ(1)).$$
Since $\NS(Z)$ is contained in the (0,0)-part
in the Hodge decomposition of $H^2(Z, \QQ(1))$,
the transcendental part $V(Z)$ inherits
a rational Hodge structure from $H^2(Z, \QQ(1))$.
Moreover,
$V(Z)$ is an irreducible Hodge structure
(\cite{Zarhin_Hodge}, Theorem 1.4.1).
Let $$\mathcal{E} = \End_{\Hdg} V(Z)$$
be the set of linear endomorphisms of $V(Z)$
that respect the Hodge structure.
Then
$\mathcal{E}$ is either a totally real field
or a CM field
(op.cit., Theorems 1.5.1 and 1.6).
With respect to the structure morphism
$Z \to {\rm Spec}\ \CC$,
there is a natural embedding
of $\mathcal{E}$ into $\CC$
given by
\begin{equation}\label{EtoC}
\mathcal{E} \to \End (V(Z) \tensor \CC)
	\to \End H^0(Z, \Omega^2) = \CC,
\end{equation}
where the second map
is the projection via the Hodge decomposition.\medskip

{\bf Corollary 2.5.}
{\em Let $k=\FF_q$ be a finite field of characteristic $p \geq 5$.
Let $X$ be a K3 surface of finite height over $k$.
Then for any quasi-canonical lifting $Z$ over $\CC$ of $X$,
the endomorphism algebra ${\cal E}$
of the transcendental cycles $V(Z)$, as a Hodge module, is a CM field
over $\QQ$.}\medskip

\pf
By Theorem 2.2,
the Tate conjecture is true for $X$ over $k$.
Since $\sigma \in \End V(Z)$
is a Hodge cycle,
$\QQ[\sigma] \rc \mathcal{E}$.
Notice that
the element $\sigma$ can not be totally real,
for otherwise
a characteristic root of $\sigma$ on $V(Z)$
would be a root of unity. This
contradicts the condition (2) in Definition 2.1.
Thus $\mathcal{E}$ must be a CM field.
\qed

\section{The Frobeius on the transcendental part}

{\bf Definition 3.1}
(cf. \cite{Zarhin_TC}, 2.0.1){\bf .}
Let $X$ be a K3 surface over $k=\FF_q$.
Let $A_{\ell}(X)$
be the set of elements
$\alpha \in H^2_{et}(X_{\bar{k}}, \QQ_{\ell}(1))$
such that $\alpha$ is invariant
under $\Gal(\bar{k}/k')$
for some finite extension $k'$ of $k$
(depending on $\alpha$).
We call $A_{\ell}(X)$
the {\em algebraic part}
of the \'etale cohomology
$H^2_{et}(X_{\bar{k}}, \QQ_{\ell}(1))$
of $X_{\bar{k}}$.
Let $V_{\ell}(X)$ be the orthogonal complement
of $A_{\ell}(X)$
with respect to the cup-product.
We call $V_{\ell}(X)$
the {\em transcendental part}
of $H^2_{et}(X_{\bar{k}}, \QQ_{\ell}(1))$.
(Cf. \cite{Zarhin_TC}, Remark 3.3.3.)
We have the decomposition
$$H^2_{et}(X_{\bar{k}},\QQ_{\ell}(1))=V_{\ell}(X)\oplus A_{\ell}(X).$$
We will write $R(X,T)$
for the reciprocal characteristic polynomial
of the geometric Frobenius $\sigma$
on $V_{\ell}(X)$.
That is,
\[ R(X,T) = \det\left(1-T\sigma | V_{\ell}(X) \right). \]\medskip

If we let
$F(X,T)$ be the reciprocal characteristic polynomial
of $\sigma$ on $H^2_{et}(X_{\bar{k}}, \QQ_{\ell}(1))$,
and put
\[ R'(X,T) = \frac{F(X,T)}{R(X,T)}. \]
Then the zeros of $R'(X,T)$
consists of zeros of $F(X,T)$
which are roots of unity.
Moreover, we have
$$R(X,qT) \in 1 + T\ZZ[T].$$
This is because all reciprocal roots of $R(X,qT)$
are algebraic integers
(cf. \cite{Zarhin_TC}, Remark 2.0.3, 2.1).
If the Tate conjecture is true for $X$ over $k$,
then $R(X,T)$
is of degree $22 - \rho(X_{\bar{k}})$,
which is the case when $k$ is of characteristic $p \geq 5$
and $X$ is of finite height
(see (2.2)).\medskip

{\bf Proposition 3.2.}
{\em Let $k=\FF_q$ be a finite field of characteristic $p$.
Let $X$ be a K3 surface of finite height $h$ over $k$.
Then $R(X,T)$ is a power of a $\QQ$-irreducible polynomial
of the form:
$$
R(X,T)=Q(X,T)^r
$$
where $Q(X,T)$ is a $\QQ$-irreducible polynomial with constant term $1$,
and the exponent $r$ divides $h$.}\medskip

\pf
We give two proofs here.
The first method uses the formal Brauer group $\hBr(X)$ associated to $X$,
and the second method is by lifting $X$ to characteristic zero
when $p \geq 5$.

(a)
Suppose that $$R(X,T) = Q_1(T) \cdots Q_r(T)$$
is the decomposition of $R(X,T)$
into irreducible polynomials with contant terms $1$ over $\QQ$.
Then over $\QQ_p$,
each $Q_i(T)$ decomposes
as a product
\[ Q_i^{<0}(T) Q_i^0(T) Q_i^{>0}(T) \]
according to the slopes of the Newton polygon of $Q_i(T)$, e.g.,
$Q_i^{<0}(T)$ denotes the polynomial corresponding to the Newton slope $<0$ part,
and respectively for $Q_i^0(T)$ and $Q_i^{>0}(T)$.
By the symmetry of the Newton polygon,
$$Q_i^{>0}(T) = (-1)^{d_i} c_i T^{d_i} Q_i^{<0} (1/T),$$
where $d_i$ is the degree of $Q_i^{<0}(T)$
and $c_i$ is the product of roots of $Q_i^{<0}(T)$.
Therefore, for each $i$,
the factor $Q_i^{<0}(T) \neq 1$,
for otherwise $Q_i(T)=Q_i^0(T)$ and
roots of $Q_i(T) = Q_i^0(T)$ are roots of unity.

On the other hand,
the product
\[ Q_1^{<0}(qT) \cdots Q_r^{<0}(qT) \]
is the reciprocal characteristic polynomial
of the Frobenius endomorphism
on the Cartier module
of $\hBr(X)$,
which is equal to $(Q^{<0}(qT))^r$
for some irreducible polynomial $Q^{<0}(qT)$
over $\ZZ_p$
(see \cite{Haz}, Theorem (24.2.6)).
Thus for all $i$, $Q_i(T)$ must be equal to each other.

(b)
Assume that $k$ is of characteristic $p \geq 5$.
Let $Z$ be a quasi-canonical lifting of $X$
to $\CC$ and $V(Z) \rc H^2(Z, \QQ(1))$ be the transcendental cycles of $Z$.
Then the Frobenius endomorphism $\sigma$ can be regarded as a Hodge cycle
in the endomorphism algebra $\mathcal{E}$
of the Hodge structure $V(Z)$.
Since $\mathcal{E}$ is a field,
the minimal polynomial $m(T)$ of $\sigma$ on $V(Z)$
has only simple roots.
Thus $m(T)$ is irreducible over $\QQ$.
Moreover,
since $V(Z)$ is also a vector space
over $\QQ[\sigma] \isom \QQ[T]/m(T)$,
the characteristic polynomial
of $\sigma$ on the $\QQ$-vector space $V(Z)$
is then the $r$-th power of $m(T)$
where $r = \dim_{\QQ[\sigma]} V(Z)$.
Therefore the assertion follows.
\qed

{\bf (3.3)}
Suppose that $X$ is of finite height $h$.  Let
$\tau = \dim_{\QQ_{\ell}} V_{\ell}(X)$.
Notice that
the splitting field
of the polynomial $Q(X,T)$
of degree $(\tau/r)$
is a CM-field
(see the proof of Corollary 2.5).
Thus $\tau/r$ is an even integer.
All possible values for $r$
in various situations
are tabulated below.
Notice that
when $X$ is ordinary (i.e., $h = 1$),
the polynomial $R(X,T)$
is always irreducible over $\QQ$.
This is a result of Zarhin
(\cite{Zarhin_TC}, Theorem 1.1)
and is used in his proof
of the Tate conjecture
for powers of an ordinary K3 surface over a finite field.\\

\vspace{.5cm}
\begin{tabular}{|c|cccccccccc|}
\hline
$\tau\backslash h$ & 1 & 2 & 3 & 4 & 5 & 6 & 7 & 8 & 9 & 10 \\
\hline
2 & 1 &&&&&&&&& \\
\hline
4 & 1 & 1, 2 &&&&&&&& \\
\hline
6 & 1 & 1 & 1, 3 &&&&&&&\\
\hline
8 & 1 & 1, 2 & 1 & 1, 2, 4 &&&&&&\\
\hline
10 & 1 & 1 & 1 & 1 & 1, 5 &&&&&\\
\hline
12 & 1 & 1, 2 & 1, 3 & 1, 2 & 1 &
	$\begin{array}{c} 1, 2, \\ 3, 6\end{array}$ &&&&\\
\hline
14 & 1 & 1 & 1 & 1 & 1 & 1 & 1, 7 &&&\\
\hline
16 & 1 & 1, 2 & 1 & 1, 2, 4 & 1 & 1, 2 & 1 &
	$\begin{array}{c} 1, 2,\\ 4, 8 \end{array}$ &&\\
\hline
18 & 1 & 1 & 1, 3 & 1 & 1 & 1, 3 & 1 & 1 & 1, 3, 9 &\\
\hline
20 & 1 & 1, 2 & 1 & 1, 2 & 1, 5 & 1, 2 & 1 & 1, 2 & 1 &
	$\begin{array}{c} 1, 2, \\ 5, 10 \end{array}$\\
\hline
\multicolumn{11}{c}{}\\
\end{tabular}\\

As an application of Proposition 3.2,
we have the following consequence,
which is communicated to us by Y.G. Zarhin.\medskip

{\bf Corollary 3.4.}
{\em Let $X_1$ and $X_2$ be two K3 surfaces
of finite height over $k = \FF_q$.
Suppose that the heights of $X_i$ are different.
Then the Tate conjecture is true for the product
$X_1 \times X_2$ over $k$.}\medskip

\pf
It suffices to show that
the tensor product $V_{\ell}(X_1) \tensor V_{\ell}(X_2)$
of the transcendental parts $V_{\ell}(X_i)$ of $X_i$
does not create non-trivial Tate classes.
For $i = 1, 2$,
let $R(X_i,T) = Q(X_i,T)^{r_i}$,
where $Q(X_i,T)$ are $\QQ$-irreducible polynomials with constant terms $1$.
Notice that $Q(X_1,T)$ and $Q(X_2,T)$ have different Newton polygons
by the assumption.
In particular $Q(X_1,T) \neq Q(X_2,T)$.

Suppose that there is a non-trivial Tate class
in $V_{\ell}(X_1) \tensor V_{\ell}(X_2)$.
Then by replacing $k$ by a finite extension,
we may assume that there exists a root $\alpha_i$ of $Q(X_i,T)$ ($i=1,2$)
such that $\alpha_1 \alpha_2 = 1$.
Thus $\alpha_2 = \alpha_1^{-1}$.
Since all the roots of $Q(X_2,T) \in \QQ[T]$ have complex norm one,
$Q(X_2, \alpha_1^{-1}) = 0$ implies that $Q(X_2, \alpha_1) = 0$.
So $Q(X_1,T)$ and $Q(X_2,T)$ have a common root.
Now the irreducibility of both polynomials $Q(X_1,T)$ and $Q(X_2,T)$
with the same constant term $1$ implies that $Q(X_1,T)=Q(X_2,T)$,
which is a contradiction.
\qed

{\bf Lemma 3.5.}
{\em Let $k=\FF_q$ be a finite field of characteristic $p \geq 5$.
Let $X$ be a K3 surface of finite height over $k$.
Suppose that $R(X,T)$ is irreducible over $\QQ$.
Then the Tate conjecture is true for the self-product
$X \times X$ of $X$ over $k$.}\medskip

\pf
This is a special case of \cite{Zarhin_power}, Theorem 4.4.
(See also op.cit., Lemma 2.1.)
Since $R(X,T)$ is irreducible,
it implies that the vector space of Galois invariants
\[ \left(V_{\ell}(X) \tensor V_{\ell}(X)\right)^{\Gal(\bar{k}/k)}
 \rc V_{\ell}(X) \tensor V_{\ell}(X) \]
is generated by the linearly independent set
$\{id, \sigma, \dots, \sigma^{\tau-1} \}$,
where $\tau = 22 - \rho(X_{\bar{k}})$.
\qed

{\bf Corollary 3.6.}
{\em Suppose that $X$ is a K3 surface of finite height
defined over the prime field $k=\FF_p$
with $p \geq 5$.
Then the Tate conjecture is true
for $X \times X$ over $k$.}\medskip

\pf
It suffices to show that $R(X,T)$ is irreducible over $\QQ$.
Notice that $R(X,T)$ is irreducible
if the action of the Frobenius $\sigma$ on the Cartier module
of $\hBr(X)$ is irreducible
(see the proof (a) of Proposition 3.2).
Under the assumption,
$\QQ_p[\sigma]$,
regarded as in $(\End_{\bar{k}} \hBr(X)) \tensor \QQ_p$,
is a totally ramified field extension
of degree $h$ over $\QQ_p$,
where $h$ is the height of $\hBr(X)$
(\cite{Haz}, Remark (24.2.7)).
Thus $R(X,T)$ is irreducible over $\QQ$.
The assertion now follows from the lemma above.
\qed

{\bf Remark 3.7.}
Let $X$ be a K3 surface of finite height over $k = \FF_q$.
As in Lemma 3.5,
we assume that $R(X,T)$ is irreducible.
If the Zariski closure of the cyclic group
generated by $\sigma/q$ is equal to
the group, $U$, of all elements of norm 1 in $\QQ[\sigma]$,
then we can use the same argument
as in \cite{Zarhin_power}, Theorem 6.1 (ii)
to show that the Tate conjecture holds true
for the $n$-th power $X^n:=X \times X \times \cdots \times X$
of $X$ over $k$ for any positive integer $n$.
However, we are not able to determine
if $\sigma/q$ generates $U$ or not.
We thank Zarhin
for pointing out to us
that the irreducibility of $R(X,T)$
is not enough to guarantee the Tate conjecture
for the $n$-th power $X^n$ of $X$ when $n>2$.

\section{Formal group laws}

As the formal Brauer group $\hBr(X)$
arising from a K3 surface $X$ over a finite field
is one-dimensional,
it must be associated to a certain Dirichlet series (\cite{Honda}).
We will discuss $\hBr(X)$
from the point of view of formal group {\it laws},
along the line of Honda \cite{Honda}.\medskip

{\bf (4.1)}
Let $k=\FF_q$ be a finite field of characteristic $p$. Let $W(k)$
be the ring of Witt vectors of $k$. Let
$X$ be a K3 surface over $k$.
Let $P(X,T)$ be the reciprocal characteristic polynomial
of the Frobenius endomorphism
on the crystalline cohomology $H^2_{cris} (X/W)$.
Then $P(X,T) \in 1 + T \cdot \ZZ [T]$,
and over $\ZZ_p$,
the polynomial $P(X,T)$
has a natural slope decomposition
\[ P(X,T) = P_{<1}(X,T) \cdot P_1(X,T) \cdot P_{>1}(X,T), \]
where $P_{<1}$ denotes the polynomial corresponding to
the Newton slope less than $1$ part,
and respectively for $P_1$ and $P_{>1}$.
If $X$ is of finite height $h$,
then $P_{<1}$ is of degree $h$
and
\[ P_{>1}(X,T) = cT^h \cdot P_{<1}(X,q^2/T) \]
for some constant $c \in \ZZ_p$.
If $X$ is supersingular,
$P_{<1} = P_{>1} = 1$.
\medskip

{\bf Theorem 4.2.}
{\em Let $X$ be a K3 surface over $k=\FF_q$.
Let
$$
P(X,T) = P_{<1}(X,T) \cdot P_1(X,T) \cdot P_{>1}(X,T)
$$
be the slope decomposition over $\ZZ_p$
of the reciprocal characteristic polynomial
of the geometric Frobenius
on $H^2_{cris}(X/W)$.
Then the isomorphism class of $\hBr(X)$ over $k$
is determined by $P_{<1}(X,T)$.}\medskip

\pf
Assume that
$X$ is of height $h < \infty$.
Let $R_{<1}(T)$
be the minimal reciprocal polynomial
of the Frobenius endomorphism
of $\hBr(X)$ on its Cartier module.
Then
$R_{<1}(T)$ is irreducible over $\ZZ_p$
(see \cite{Haz}, Theorem (24.2.6))
and
\[ P_{<1}(X,T) = R_{<1}(T)^r \]
for some positive integer $r$.
One knows that
the formal group $\hBr(X)$
is determined by $R_{<1}(T)$
(op.cit., Proposition (24.2.9)),
and thus it is determined by $P_{<1}(X,T)$.

On the other hand,
if $X$ is supersingular,
then $P_{<1} = 1$
and $\hBr(X)$ is isomorphic to
the additive formal group $\hat{\GG}_a$.
The assertion is trivial in this case.
\qed

{\bf (4.3)}
In what follows, we let $k = \FF_q$ be the finite field of $q$ elements
of characteristic $p$, and
let $W = W(k)$ be the ring of Witt vectors of $k$.\medskip

Here we compute explicitly the formal group laws
that realize the formal Brauer groups $\hBr(X)$ of certain K3 surfaces $X$
over $k$.
The explicit realizations are obtained
based on the work of Stienstra \cite{St}.\medskip

The examples below involve hypergeometric series.
This seems to come from the fact that
the principle part of the defining equation
for a K3 surface discussed in our examples is very
symmetric with respect to the variables,
and the parameter $t$
appears in the product of the variables.
We recall that
the hypergeometric series ${}_mF_n$
with upper parameters $\{a_i\}_{i=1}^m$
and lower parameters $\{b_i\}_{i=1}^n$,
$a_i, b_i \in \CC$ and $b_i \not\in \ZZ_{\leq 0}$,
is defined by the formal power series
\[ {}_mF_n\left(\begin{array}{c}
	a_1, a_2, \cdots, a_m \\ b_1, b_2, \cdots, b_n
	\end{array}; x \right)
	= \sum_{r=0}^{\infty}
	\frac{(a_1)_r (a_2)_r \cdots (a_m)_r}{(b_1)_r (b_2)_r \cdots (b_n)_r}
	\cdot \frac{x^r}{r!}, \]
where $(a)_0 = 1$ and
$(a)_r = a(a+1) \cdots (a+r-1)$ for $r > 0$
is the {\em Pochhammer symbol}.
Notice that
if $a_i$ is a non-positive integer for some $i$,
then the hypergeometric series ${}_mF_n$
is a polynomial.
We may regard that the series
is defined over a ring $R$
whenever the expansion makes sense over $R$.\\

{\bf (4.4) Elliptic modular K3 surfaces.}\\

{\bf (4.4.1)} (Shioda \cite{Shioda})
Let
\[ \Gamma(4)= \left\{\left(\begin{array}{cc} a & b \\ c & d \end{array}\right)
\mid c\equiv 0 \pmod 4\right\} \]
be the principal congruence subgroup
of $\Gamma=\text{SL}_2(\ZZ)/\{\pm 1\}$ of level $4$.
Then $\Gamma(4)$ gives rise to an elliptic modular K3 surface
$X$ over a finite field $k=\FF_q$ of characteristic $p>2$
with $q = p^a$.
More precisely,
$X$ is the minimal model over $k$ of the Jacobi quartic
$$
y^2=(1-\sigma^2x^2)(1-\sigma^{-2}x^2),
$$
which is defined over the function field $K=k(\sigma)$
with $\sigma$ a variable over $k$.
Suppose $\sqrt{-1} \in k$.
Then the zeta-function of $X$ over $k$
is given as follows:
$$
Z(X,T)=\frac{1}{(1-T)(1-qT)^{20}(1-q^2T)}
	\cdot \frac{1}{H_{3,q}(T)}.
$$
Here $H_{3,q}(T)$ is the Hecke polynomial
$$
H_{3,q}(T)= \left\{
\begin{array}{ll}
	1-(\pi^2+\pi'^2)T+q^2T^2 & \text{if}\quad
		p\equiv 1\pmod 4\\
	(1-qT)^2 & \text{if}\quad p\equiv 3\pmod 4
\end{array} \right.
$$
where $\pi, \pi'$ are integers in $\ZZ[\sqrt{-1}]$ such that
$\pi = \alpha^a$ and $\pi = \alpha'^a$
for some $\alpha, \alpha' \in \ZZ[\sqrt{-1}]$
with $\alpha \alpha' = p$
and $\alpha \equiv 1\pmod{2\sqrt{-1}}$
(\cite{Shioda}, Example A.10).
Fix an embedding
$\bar{\QQ} \hookrightarrow \bar{\QQ}_p$.
When $p \equiv 1 \pmod{4}$, assume that
$\pi$ is a $p$-adic unit
with respect to this embedding.\medskip

{\bf (4.4.2)}
Make the change of variable
$x \mapsto \sigma x$ in the defining equation for the Jacobi quartic. Then
this leads to the elliptic pencil given by
\[ y^2 = (1 - \sigma^4 x^2)(1-x^2) \]
with parameter $\sigma$.
Consider the formal group law $\hat{G}$
over $\ZZ$
with the logarithm
\[ l(\tau) = \sum_{m=0}^{\infty} a(m) \frac{\tau^{4m+1}}{4m+1}, \]
where
\begin{eqnarray*}
a(m) &=& \binom{2m}{m} \cdot {}_2F_1\left(\begin{array}{c}
						-m, -m \\ 1 \end{array}; 1\right) \\
	&=& \binom{2m}{m}^2.
\end{eqnarray*}
The last equality follows
from the Vandermonde involution formula.
Then the base change to $k$ of $\hat{G}$
realizes $\hBr(X)$
(\cite{St}, Theorem 2).\medskip

Suppose $p = 4m+1$.
Then $X$ is of height one.
Take $k = \FF_p$ to be the prime field.
Then for any integer $s \geq 0$,
we have (\cite{SB}, Theorem (A.8)(v))
\[ \frac{a\left(\frac{p^{s+1}-1}{4}\right)}{a\left(\frac{p^s-1}{4}\right)}
	\equiv \pi^2 \pmod{p^{s+1}}, \]
where $\pi^2$ is the $p$-adic unit reciprocal root
of $H_{3,q}(T), q=p$ discussed in (4.4.1).
On the other hand,
since the polynomials
\[ {}_2F_1\left(\begin{array}{c}
		-m_s, -m_s \\ 1 \end{array}; 1\right), \]
with $m_s: = \frac{1}{4}(p^s - 1)$, converge $p$-adically to
\[ {}_2F_1\left(\begin{array}{c}
		\frac{1}{4}, \frac{1}{4} \\ 1 \end{array}; 1\right), \]
one knows that
\[ \lim_{s \to \infty}
	\frac{a\left(\frac{p^{s+1}-1}{4}\right)}{a\left(\frac{p^s-1}{4}\right)}
	= h(1)^2. \]
Here $h(x)$ is the formal power series
\[ h(x) = {}_2F_1\left(\begin{array}{c}
			\frac{1}{4}, \frac{1}{4} \\ 1 \end{array}; x\right)
	\Big/ {}_2F_1\left(\begin{array}{c}
			\frac{1}{4}, \frac{1}{4} \\ 1 \end{array}; x^p\right) \]
with coefficients in $\ZZ_p$, and it converges $p$-adically at $x = 1$.
Furthermore,
since
\[ -\frac{1}{4} = \sum_{i=0}^{\infty} mp^i, \]
we have
\[ h(1) = \frac{\Gamma_p(\frac{1}{4})^2}{\Gamma_p(\frac{1}{2})}, \]
where $\Gamma_p(x)$
is the Morita's $p$-adic gamma-function
(\cite{Koblitz_Hyper}, Theorem 2).
Thus
\[ \pi^2 = \frac{\Gamma_p(\frac{1}{4})^4}{\Gamma_p(\frac{1}{2})^2}. \]
This is a consequence of the Gross-Koblitz formula
(see \cite{vH}, (F.2)
combined with some basic equalities
for $\Gamma_p(x)$ in \cite{Koblitz_p}, p.42).\medskip

If $p \equiv -1 \pmod{4}$,
then $X$ is supersingular.
Thus $l(\tau) \in W[[\tau]]$,
that is,
\[ p^{2s}\ \bigg|\ a\left(\frac{p^{2s}-1}{4}\right), \]
which implies that
\[ p^s\ \bigg|\ \binom{\frac{p^{2s}-1}{2}}{\frac{p^{2s}-1}{4}} \]
for all non-negative integers $s$.\\

{\bf (4.5) Examples: Diagonal quartic K3 surfaces}\\

{\bf (4.5.1)}
Consider the twisted diagonal quartic surface
$X$ over $k$
defined by
$$
c_1T_1^4+c_2T_2^4+c_3T_3^4+c_4T_4^4=0 \rc
	\mathbb{P}^3_k
$$
with $c_i \in k$ and $c = c_1c_2c_3c_4\ne 0$.
Let $\hat{c_i} \in W$
be a lifting of $c_i$ for each $i$
and let
\[ \hat{c} = \hat{c_1} \hat{c_2} \hat{c_3} \hat{c_4}. \]
Put
$$
l(\tau) =\sum_{n=0}^{\infty}
	\hat{c}^n\frac{(4n)!}{(n!)^4} \frac{\tau^{4n+1}}{4n+1}
$$
and define
$$
\hat{G}(\tau_1,\tau_2)=l^{-1}(l(\tau_1)+l(\tau_2)).
$$
Then $\hat{G}$ is a formal group
(defined over $W$)
whose reduction $G$ to $k$
realizes $\hBr(X)$
(\cite{St}, Theorem 1).\medskip

Now let $k=\FF_q$ be the finite field
with $q=p^a$ elements of characteristic $p$.
The formal group $G(\tau_1,\tau_2)$
is of multiplicative type (i.e., $G$ is of height 1,
and hence $G$ is isomorphic over $\bar{k}$ to $\hat\GG_m$)
if $p\equiv 1 \pmod 4$, and is additive (i.e., $G$ has infinite height and
hence $G$ is isomorphic to $\hat\GG_a$) if $p\equiv 3\pmod 4$.
In the later case,
it means that
$l(\tau)$ is an element in $W[[\tau]]$.
Thus it implies that
for all $n > 1$,
we have
\[ \frac{(4n)!}{(n!)^4} \equiv 0 \pmod{p^{s}}, \]
where $s$ is the largest integer such that $p^s\ |\ (4n+1)$.\medskip

On the other hand,
assume that $p \equiv 1 \pmod{4}$.
Write
\[ a(n) = \hat{c}^n \frac{(4n)!}{(n!)^4}. \]
Then there exists an element $\alpha \in W$
such that
\begin{equation}\label{a1}
a\left( \frac{\mu p^{s+1} -1}{4} \right) \equiv
	\alpha \cdot a\left( \frac{\mu p^s -1}{4}\right) \pmod{p^{s+1}}
\end{equation}
for all non-negative integers $\mu, s$
(with the convention that
$a(n) = 0$ if $n$ is not an integer)
(\cite{SB}, Theorem (A.8)(v)).
Let
\[ \pi = \alpha^{1+ \sigma + \cdots +\sigma^{a-1}}. \]
Here the upper script $\sigma$
indicates the action by
the absolute Frobenius $\sigma$ on $W$.
Then $\pi$ turns out to be a twisted Jacobi sum
(\cite{GY}, Lemma 3.4)
and equation (\ref{a1}) gives
$p$-adic approximations
of $\pi$.\medskip

{\bf (4.5.2)}
Let $p$ be an odd prime.
Generalizing the above example,
we consider the one-parameter family of K3 surfaces $X_{\lam}$
given by one of the following quadrics :
\[ c_1T_1^4 + c_2T_2^4 + c_3T_3^4 + c_4T_4^4
	- 4\lam T_1T_2T_3T_4 = 0, \]
\[ c_1T_1^3T_2 + c_2T_2^3T_3 + c_3T_3^3T_4 + c_4T_4^3T_1
	- 4\lam T_1T_2T_3T_4 = 0, \]
\[ c_1T_1^2T_2T_3 + c_2T_2^2T_3T_4 + c_3T_3^2T_4T_1
	+ c_4T_4^2T_1T_2 - 4\lam T_1T_2T_3T_4 = 0, \]
with $c_i \in k$, $c = c_1c_2c_3c_4\ne 0$
and parameter $\lam$.
We only consider those values of $\lam$
such that the corresponding surfaces $X_{\lam}$ are smooth.
Let $\hat{c_i} \in W$
be a lifting of $c_i$ for each $i$,
and let
\[ \hat{c} = \hat{c_1} \hat{c_2} \hat{c_3} \hat{c_4} \in W. \]
Also take a lifting $\hat{\lam} \in W$ of $\lam$.
Put
\[ l(\tau) = \sum_{m=0}^{\infty} a(m) \frac{\tau^{m+1}}{m+1}, \]
where
\[ a(m) = (-4\hat{\lam})^m {}_4F_3\left(\begin{array}{c}
			\frac{-m}{4},\frac{-m+1}{4},\frac{-m+2}{4},\frac{-m+3}{4} \\
					1,1,1 \end{array}; \hat{c} \hat{\lam}^{-4} \right). \]
Define
\[ G_{\hat{\lam}} (\tau_1, \tau_2) = l^{-1} ( l(\tau_1) + l(\tau_2) ). \]
Then $G_{\hat{\lam}}$ is a formal group
(defined over $W$)
whose reduction $G_{\lam}$ to $k$
realizes $\hBr(X_{\lam})$
(\cite{St}, Theorem 1).\medskip

If $a(p-1) \equiv 0 \pmod{p}$,
then $G_{\lam}$ is isomorphic to the additive group $\hat{\GG}_a$.
Thus in this case, for any $m > 0$,
\[ a(m) \equiv 0 \pmod{p^s} \]
where $s$ is the largest integer such that $p^s\ |\ (m+1)$.
\medskip

On the other hand,
for those $\lam$
satisfying
$a(p-1) \not\equiv 0 \pmod{p}$,
the formal group $G_{\lam}$ is of multiplicative type, i.e., is isomorphic over
$\bar k$ to the multiplicative group $\hat\GG_m$.
In this case,
we take $\hat{c}$ and $\hat{\lam} \in W$
to be the Teichm\"uller lifting of $c$ and $\lam$,
respectively.
Then for each such $\lam$,
there exists a $p$-adic unit $\alpha_{\lam} \in W$
such that
\[ a(\mu p^{s+1}-1) \equiv
	\alpha_{\lam} \cdot a(\mu p^s - 1)^{\sigma} \pmod{p^{s+1}} \]
for any integers $\mu, s \geq 0$
(\cite{SB}, Theorem (A.8)(v)).
Furthermore,
if $q = p^a$,
the element
\[ \pi_{\lam} = \alpha_{\lam}^{1+\sigma+ \cdots + \sigma^{a-1}} \]
is the unique $p$-adic unit root
of the characteristic polynomial
of the geometric Frobenius
on $H^2_{et} ((X_{\lam})_{\bar{k}}, \QQ_{\ell})$.
It is then easy to see that
the polynomials
\[ {}_4F_3\left(\begin{array}{c}
	\frac{-m_s}{4},\frac{-m_s+1}{4},\frac{-m_s+2}{4},\frac{-m_s+3}{4} \\
		1,1,1 \end{array}; x \right) \]
with $m_s: = (p^s - 1)$, converge $p$-adically to
the hypergeometric series
\[ {}_3F_2 \left(\begin{array}{c}
	\frac{1}{4}, \frac{2}{4}, \frac{3}{4} \\
	1,1 \end{array}; x \right)
	= \sum_{r=0}^{\infty} 4^{-4r}
		\binom{4r}{r}\binom{3r}{r}\binom{2r}{r} x^r. \]
(Notice that
this formal power series ${}_3F_2$
has coefficients in $\ZZ[\frac{1}{2}]$.)
Therefore
\begin{eqnarray*}
\alpha_{\lam} &=& \lim_{s \to \infty} \frac{a(p^{s+1}-1)}{a(p^s-1)^{\sigma}} \\
	&=& \hat{\lam}^{p-1} f(\hat{c}\hat{\lam}^{-4}).
\end{eqnarray*}
Here $f(x)$ is the formal power series
\begin{equation}\label{f1}
f(x) = {}_3F_2 \left(\begin{array}{c}
	\frac{1}{4}, \frac{2}{4}, \frac{3}{4} \\
	1,1 \end{array}; x \right) \Big/
	{}_3F_2 \left(\begin{array}{c}
	\frac{1}{4}, \frac{2}{4}, \frac{3}{4} \\
	1,1 \end{array}; x^p \right)
\end{equation}
with coefficients in $\ZZ_p$, and it
converges $p$-adically at $x=\hat{c}\hat{\lam}^{-4}$.
Note that
the requirement of $\hat{c}\hat{\lam}^{-4}$
being the Teichm\"uller lifting of $c\lam^{-4}$
is needed
in order to have the compact formula (\ref{f1})
for $f(x)$.\\

{\bf (4.6) Examples: Double sextic K3 surfaces}\\

{\bf (4.6.1)}
Let $X$ be the double cover of a smooth sextic
defined in $\PP^1 \times \PP^2$ over $k$ by
\[ Y^2 = c_1T_1^6 + c_2T_2^6 + c_3T_3^6 \]
with $c_i \in k$ and $c = c_1c_2c_3 \ne 0$.
Let $\hat{c_i} \in W$
be a lifting of $c_i$ for each $i$
and let
\[ \hat{c} = \hat{c_1} \hat{c_2} \hat{c_3}, \]
which is a lifting of $c$ to $W$.
Put
$$
l(\tau)=\sum_{n=0}^{\infty}
	\hat{c}^n\frac{(3n)!}{(n!)^3} \frac{\tau^{6n+1}}{6n+1}
$$
and
$$
\hat{G}(\tau_1,\tau_2) = l^{-1}(l(\tau_1)+l(\tau_2)).
$$
Then $\hat{G}$ is a formal group
(defined over $W$)
whose reduction $G$ to $k$
realizes $\hBr(X)$
(\cite{St}, Theorem 2).\medskip

Over the finite field $k$ of characteristic $p$,
the formal group $G(\tau_1,\tau_2)$ is isomorphic over $\bar{k}$
to the multiplicative group $\hat \GG_m$ if $p\equiv 1 \pmod 6$,
and to the additive group $\hat \GG_a$ otherwise.
In the later case, it follows that
$l(\tau)$ is an element in $W[[\tau]]$.
Thus it implies that for all $n > 1$,
we have
\[ \frac{(3n)!}{(n!)^3} \equiv 0 \pmod{p^{s}} \]
where $s$ is the largest integer such that $p^s\ |\ (6n+1)$.\medskip

{\bf (4.6.2)}
More generally, let $X_{\lam}, \lam \neq 0$
be the one-parameter family of K3 surfaces
given by the double cover of $\mathbb{P}^2$ ramified along
any one of the following curves :
\[ c_1T_1^6 + c_2T_2^6 + c_3T_3^6
	- 3\lam T_1^2 T_2^2 T_3^2 =0, \]
\[ c_1T_1^5T_2 + c_2T_2^5T_3 + c_3T_3^5T_1
	- 3\lam T_1^2 T_2^2 T_3^2= 0, \]
\[ c_1T_1^4T_2^2 + c_2T_2^4T_3^2 + c_3T_3^4T_1^2
	- 3\lam T_1^2 T_2^2 T_3^2= 0, \]
\[ c_1T_1^3T_2^2T_3 + c_2T_2^3T_3^2T_1 + c_3T_3^3T_1^2T_2
	- 3\lam T_1^2 T_2^2 T_3^2 = 0, \]
with $c_i \in k$, $c = c_1c_2c_3 \ne 0$
and parameter $\lam$.
We suppose that $k$ is of characteristic $p > 3$.
We only consider those values of $\lam$
such that the corresponding surfaces $X_{\lam}$ are smooth.
Let $\hat{c_i} \in W$
be a lifting of $c_i$ for each $i$
and let
\[ \hat{c} = \hat{c_1} \hat{c_2} \hat{c_3}. \]
Also take a lifting $\hat{\lam} \in W$ of $\lam$.
Put
\[ l(\tau) = \sum_{m=0}^{\infty} b(m)
	\frac{\tau^{2m+1}}{2m+1}, \]
where
\[ b(m) = (-3\hat{\lam})^m
	{}_3F_2\left(\begin{array}{c}
		\frac{-m}{3},\frac{-m+1}{3},\frac{-m+2}{3} \\
		1,1 \end{array}; \hat{c}\hat{\lam}^{-3} \right). \]
Define
\[ G_{\hat{\lam}} (\tau_1, \tau_2) = l^{-1} ( l(\tau_1) + l(\tau_2) ). \]
Then $G_{\hat{\lam}}$ is a formal group
(defined over $W$)
whose reduction $G_{\lam}$ to $k$
realizes $\hBr(X_{\lam})$
(\cite{St}, Theorem 2).\medskip

We have the similar consequences
as in (4.5.2).
Suppose that
$b(p-1) \not\equiv 0 \pmod{p}$
for some $\lam$
(i.e., $G_{\lam}$ is of multiplicative type).
Take $\hat{c}, \hat{\lam} \in W$
to be the Teichm\"uller lifting of $c, \lam$,
respectively.
Let $\beta_{\lam} \in W$
be the $p$-adic unit
such that
\[ b\left(\frac{\mu p^{s+1}-1}{2}\right) \equiv
	\beta_{\lam} \cdot b\left(\frac{\mu p^s - 1}{2}\right)^{\sigma}
		\pmod{p^{s+1}} \]
for any integers $\mu, s \geq 0$
(with the convention that
$b(m) = 0$ if $m$ is not an integer).
Here the upper script $\sigma$
indicates the action by
the absolute Frobenius $\sigma$ on $W$.
In this case, the polynomials
\[ {}_3F_2\left(\begin{array}{c}
		\frac{-m_s}{3},\frac{-m_s+1}{3},\frac{-m_s+2}{3} \\
		1,1 \end{array}; x \right) \]
with $m_s: = \frac{1}{2}(p^s - 1)$, converge $p$-adically to
\[ {}_3F_2 \left(\begin{array}{c}
	\frac{1}{6}, \frac{3}{6}, \frac{5}{6} \\
	1,1 \end{array}; x \right)
	= \sum_{r=0}^{\infty} 12^{-3r}
		\binom{6r}{r}\binom{5r}{r}\binom{4r}{r} x^r. \]
(Notice that
this formal power series ${}_3F_2$
has coefficients in $\ZZ[\frac{1}{6}]$.)
Thus we have
\begin{eqnarray*}
\beta_{\lam} &=& \lim_{s \to \infty}
	\frac{a\left(\frac{p^{s+1}-1}{2}\right)}
		{a\left(\frac{p^s-1}{2}\right)^{\sigma}} \\
	&=& (-3\hat{\lam})^{(p-1)/2} \cdot g(\hat{c}\hat{\lam}^{-3}).
\end{eqnarray*}
Here $g(x)$ is the formal power series
\[ g(x) = {}_3F_2 \left(\begin{array}{c}
	\frac{1}{6}, \frac{3}{6}, \frac{5}{6} \\
	1,1 \end{array}; x \right) \Big/
	{}_3F_2 \left(\begin{array}{c}
	\frac{1}{6}, \frac{3}{6}, \frac{5}{6} \\
	1,1 \end{array}; x^p \right) \]
with coefficients in $\ZZ_p$, and it
converges $p$-adically at $\hat{c}\hat{\lam}^{-3}$.\medskip

{\bf Remark 4.7.}
Note that
if $c_i \in \QQ$
in Examples (4.5.1) and (4.6.1),
the K3 surfaces $X$ over $\QQ$
have geometric Picard rank $\rho(X_{\bar{\QQ}}) = 20$.
Thus there is a modular form associated to
the Galois representation
on the ($\ell$-adic) transcendental cycles of each such $X$.
For this modularity property,
see \cite{Y},
especially \S 4 and Proposition 8.12.\\

{\bf (4.8) Examples of K3 surfaces of finite height $>1$}\\

Finally we give two examples
of K3 surfaces of finite height greater than 1.
We introduce some notations.
For a positive integer $n$,
denote
\[ [n] = \left\{\frac{i}{n}\ \big|\ 1 \leq i \leq n, (i,n) = 1 \right\}. \]
If $m$ is another positive integer,
$[n]^m$ denotes the multi-index
such that each element in $[n]$ repeats $m$-times.\\

{\bf (4.8.1) Examples: Quasi-diagonal K3 surfaces}\\

Consider the family
of deformations $X_{\lam}$, $\lam \in \PP^1_k$,
of the quasi-diagonal K3 surface in $\PP^3$
\[ X_{\lam}: T_1^4 + T_1T_2^3 + T_3^4 + T_4^4
	- 12\lam T_1T_2T_3T_4= 0. \]
As before,
we only consider those values of $\lam$
such that the corresponding surfaces $X_{\lam}$ are smooth.
Then similar to (4.5.1), when $\lam = 0$,
the formal Brauer group of $X_0$
is realized as the reduction to $k$
of the formal group law $G_{\hat{0}}$ over $\ZZ$
with the logarithm
\[ l(\tau) = \sum_{n=0}^{\infty}
	\binom{12n}{2n}\binom{10n}{3n}\binom{7n}{4n}
	\frac{\tau^{12n+1}}{12n+1}. \]\medskip

Similar to (4.5.2),
when $\lam \neq 0$,
let $\hat{\lam} \in W$ be a lifting of $\lam$.
Then the formal Brauer group $G_{\lam}$ of $X_{\lam}$
is realized as the reduction to $k$
of the formal group law $G_{\hat{\lam}}$ over $W$
with the logarithm
\[ l(\tau) = \sum_{m=0}^{\infty} (-12\hat{\lam})^m
	A_m(\hat{\lam}) \frac{\tau^{m+1}}{m+1}, \]
where
\[ A_m(x) = {}_{12}F_{11}\left(\begin{array}{c}
	\left\{\frac{-m+i}{12}\ |\ 0 \leq i \leq 11\right\} \\
	{[1]}^3, [2]^2, [3]^2, [4] \end{array};
		\left(2^{10} 3^6 x^{12}\right)^{-1} \right). \]
Notice that $A_m(x)$ are polynomials over $\ZZ[\frac{1}{6}]$.\medskip

Assume that $p \neq 2, 3$.
Then for $\lam \in k$,
the formal group $G_{\lam}$ is of height one
if and only if
$\lam^{p-1}A_{p-1}(\lam) \neq 0$.
Here we regard $x^{p-1}A_{p-1}(x)$
as a polynomial in $x$ over $k$.
Thus $G_{\lam}$ is generically of height one.
Put $x = \left(2^{10} 3^6 \lam^{12}\right)^{-1}$.
Let $$V_1(x) = A_{p-1}(\lam),$$
and
\[ V_2(x) = \frac{1}{p}\left(A_{p^2-1}(\lam) - A_{p-1}(\lam)^{p+1}\right) \]
regarded as polynomials of $x$.
Then in fact $V_2(x) \in W[x]$.
Reducing modulo $p$,
we regard $V_i(x)$ as polynomials over $k$.
Then for a non-zero $\lam \in k$,
the formal group $G_{\lam}$ is of height two
if and only if
$V_1(x) = 0$ but $V_2(x) \neq 0$
(see \cite{L}, Lemma (2.1)).\medskip

For example,
when $p = 13$,
as polynomials over $k$,
we have
$$V_1(x) = 1 + 10x$$ and
\begin{eqnarray*}
V_2(x) &=& 8x(x+2)(x+6)(x+10)(x^2+8x+10) \\
	&& \hspace{.5cm}\times (x^8+9x^7+5x^6+12x^5+6x^3+3x^2+7).
\end{eqnarray*}
One checks that
$V_1$ does not divide $V_2$.
Thus for those $\lam \in k$
such that $\left(2^{10} 3^6 \lam^{12}\right)^{-1} = 9$,
the K3 surface $X_{\lam}$ is of height two.\medskip

If $m$ runs through those values $m_s: = p^s - 1$,
the polynomials $A_m$ converges $p$-adically to
\[ A(x) = {}_6F_5\left(\begin{array}{c}
	[6], [12] \\
	{[1]}^2, [2], [3] \end{array};
		\left(2^{10} 3^6 x^{12}\right)^{-1} \right). \]
Notice that $A(x)$ is a formal power series (in $x^{-1}$) over $\ZZ_p$
for all primes $p \neq 2, 3$.\\

{\bf (4.8.2) Examples: Double sextic K3 surfaces of quasi-diagonal type}\\

Consider the double cover $X_{\lam}$, $\lam \in \PP^1_k$, of $\PP^2$
ramified along
\[ T_1^6 + T_1T_2^5 + T_3^6 - 15\lam T_1^2T_2^2T_3^2=0. \]
As before,
we only consider those values of $\lam$
such that the corresponding surfaces $X_{\lam}$ are smooth.
Then similar to (4.6.1),
when $\lam = 0$,
the formal Brauer group of $X_0$
is realized as the reduction to $k$
of the formal group law $G_{\hat{0}}$ over $\ZZ$
with the logarithm
\[ l(\tau) = \sum_{n=0}^{\infty}
	\binom{15n}{5n}\binom{10n}{6n} \frac{\tau^{30n+1}}{30n+1}. \]
\medskip

Similar to (4.6.2),
when $\lam \neq 0$,
let $\hat{\lam} \in W$ be a lifting of $\lam$.
Then the formal Brauer group $G_{\lam}$ of $X_{\lam}$
is realized as the reduction to $k$
of the formal group law $G_{\hat{\lam}}$ over $W$
with logarithm
\[ l(\tau) = \sum_{m=0}^{\infty} (-15\hat{\lam})^m
	B_m(\hat{\lam}) \frac{\tau^{2m+1}}{2m+1}, \]
where
\[ B_m(x) = {}_{15}F_{14}\left(\begin{array}{c}
	\left\{\frac{-m+i}{15}\ |\ 0 \leq i \leq 14\right\} \\
	{[1]}^2, [2]^2, [3], [4], [5], [6] \end{array};
		\left(4^4 5^5 6^6 x^{15}\right)^{-1} \right). \]
Notice that $B_m(x)$ are polynomials over $\ZZ[\frac{1}{30}]$.\medskip

Assume that $p > 5$.
Let $n = \frac{1}{2}(p-1)$.
Then for $\lam \in k$,
the formal group $G_{\lam}$ is of height one
if and only if
$\lam^n B_n(\lam) \neq 0$.
Here we regard $x^n B_n(x)$
as a polynomial over $k$.
Thus $G_{\lam}$ is generically of height one.
Put $x = \left(4^4 5^5 6^6 \lam^{15}\right)^{-1}$.
Similar to the previous example,
let $$V_1(x) = B_n(\lam)$$ and
\[ V_2(x) = \frac{1}{p}\left(B_{(p^2-1)/2}(\lam) - B_n(\lam)^{p+1}\right). \]
One compute that if $p = 31$,
then over $k$
we have $$V_1(x) = 1 + 20x$$ and
\begin{eqnarray*}
V_2(x) &=& 7x+2x^2+ \cdots + 24x^{32} \\
	&=& 24x(x^3 + \cdots)(x^6 + \cdots)(x^{22} + \cdots),
\end{eqnarray*}
where the last three factors
are irreducible polynomials of degree
3, 6, and 22 over $\FF_{31}$, respectively.
We have checked that $V_1$ and $V_2$
are relatively prime to each other.
Thus for those $\lam \in k$
such that $\left(4^4 5^5 6^6 \lam^{15}\right)^{-1} = 17$,
the K3 surface $X_{\lam}$ is of height two.\medskip

For an odd prime $p$,
if $m$ runs through those values
$m_s: = \frac{1}{2} (p^s - 1)$,
the polynomials $B_m$ converges $p$-adically to
\[ B(x) = {}_{12}F_{11}\left(\begin{array}{c}
	[10], [30] \\
	{[1]}^2, [2], [3], [4], [5] \end{array};
		\left(4^4 5^5 6^6 x^{15}\right)^{-1} \right). \]
Notice that $B(x)$ is a formal power series (in $x^{-1}$)
over $\ZZ_p$ for all $p > 5$.
\medskip

{\bf Question 4.9.}
For all prime numbers $p < 150$ in the last two examples,
we have checked that
for any non-zero $\lam \in k$,
the formal group $G_{\lam}$ is either of height 1 or 2,
i.e.,
the two polynomials $V_1$ and $V_2$ over $\FF_p$
have no non-trivial common divisor.
Is it true that $G_{\lam}$ has height $\leq 2$
for each $\lam \neq 0$ and every possible prime $p$?\medskip

{\bf Acknowledgments.}
We would like to thank Y. G. Zarhin and M. Sch\"utt
for their interest in this work and for providing useful suggestions
and corrections to the early version of this paper.
Sch\"utt helped us understand the statements in (4.4.1).
Zarhin kindly pointed out to us a consequence of Proposition 3.2, which is
formulated as Corollary 3.4.

\end{document}